\documentclass[10pt]{amsart}
\usepackage{amssymb}
\usepackage{epsfig}
\usepackage{latexsym}
\usepackage{amsmath,amssymb,amsfonts,amsthm,graphics}
\usepackage{eucal}
\usepackage{eufrak}
\usepackage[all]{xypic}
\usepackage{xspace}

\theoremstyle{plain}
\newtheorem{theorem}{Theorem}

\newtheorem{lemma}{Lemma}

\theoremstyle{definition}

\theoremstyle{remark}

\numberwithin{equation}{section}
\def\Tilde{\char126\relax}

\begin{document}


\title[Efficient Counting and Asymptotics of $k$-noncrossing tangled-diagrams]
      {Efficient Counting and Asymptotics of $k$-noncrossing tangled-diagrams}
\author{William Y. C. Chen$^{\dagger}$, Jing Qin$^{\dagger}$,
        Christian M. Reidys$^{\dagger}$$^{\,\star}$ and
        Doron Zeilberger$^{\sharp}$}
\address{ $^{\dagger}$Center for Combinatorics, LPMC-TJKLC \\
          Nankai University  \\
          Tianjin 300071\\
          P.R.~China\\
          Phone: *86-22-2350-6800\\
          Fax:  *86-22-2350-9272\\
            $^{\sharp}$Mathematics Department,\\
  Rutgers University(New Brunswick),
  Piscataway, NJ, USA.\\
zeilberg@math.rutgers.edu.} \email{reidys@nankai.edu.cn}
\thanks{}
\keywords{matching, vacillating tableau, holonomic ansatz,
singularity, singular expansion, D-finite}
\date{February, 2008}
\begin{abstract}
In this paper we enumerate $k$-noncrossing tangled-diagrams. A
tangled-diagram is a labeled graph whose vertices are $1,\dots,n$
have degree $\le 2$, and are arranged in increasing order in a
horizontal line. Its arcs are drawn in the upper halfplane with a
particular notion of crossings and nestings. Our main result is the
asymptotic formula for the number of $k$-noncrossing
tangled-diagrams $T_{k}(n) \, \sim \,c_k \, n^{-((k-1)^2+(k-1)/2)}\,
(4(k-1)^2+2(k-1)+1)^n$ for some $c_k>0$.
\end{abstract}
\maketitle {{\small
}}



\section{Tangled diagrams as molecules or walks}\label{S:intro}


In this paper we show how to compute the numbers of $k$-noncrossing
tangled-diagrams and prove the asymptotic formula
\begin{equation}\label{E:formula}
T_{k}(n) \, \sim  \, c_k  \, n^{-((k-1)^2+(k-1)/2)}\,
(4(k-1)^2+2(k-1)+1)^n,\qquad c_k>0 \ .
\end{equation}
This article is accompanied by a Maple package {\tt TANGLE},
downloadable from the webpage
$$
{\tt http://www.math.rutgers.edu/ \Tilde
zeilberg/mamarim/mamarimhtml/tangled.html} \quad .
$$

\begin{figure}[ht]\label{E:III}
\begin{center}
\scalebox{0.7}[0.9]{\includegraphics*[10,740][450,810]{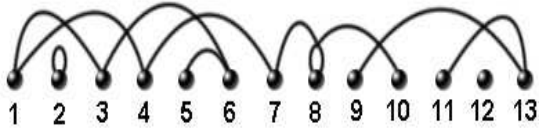}}
\end{center}
\caption{\small{A tangled-diagram over $13$ vertices.}}
\label{F:tan}
\end{figure}

$k$-noncrossing tangled-diagrams are motivated by studies of RNA
molecules. They serve as combinatorial frames for searching
molecular configurations and were recently studied
\cite{Reidys:07vac} by the first three authors. Tangled-diagrams are
labeled graphs over the vertices $1,\dots,n$, drawn in a horizontal
line in increasing order. Their arcs are drawn in the upper
halfplane having the following types of arcs
\begin{center}
\scalebox{0.6}[0.6]{\includegraphics*[60,650][560,820]{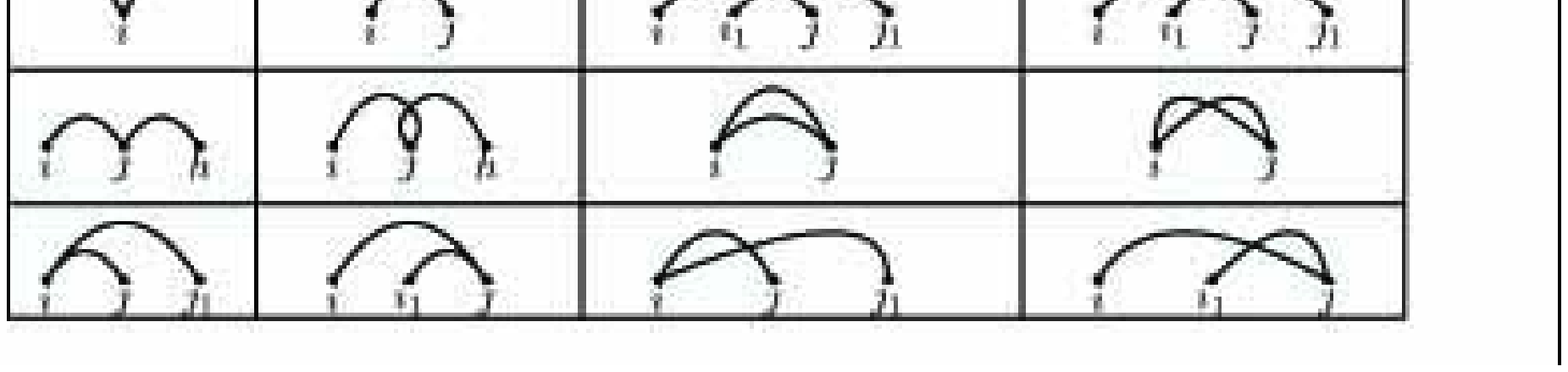}}
\end{center}
Tangled diagrams have possibly isolated points, for instance, the
tangled diagram displayed in Figure~\ref{F:tan} has the isolated
point $12$.
Details on tangled-diagrams can be found in \cite{Reidys:07vac}. It
is interesting to observe that tangled-diagrams are in
correspondence to
the following types of walks:\\

{\bf Observation 1}: The number of $k$-noncrossing tangled-diagrams
over $[n]$, without isolated points, equals the number of simple
lattice walks in $x_1 \geq x_2 \geq \dots \geq x_{k-1} \geq 0$, from
the origin back to the origin, taking $n$ days, where at each day
the walker can either make {\it one} unit step in any (legal)
direction, or else feel energetic and make any {\it two} consecutive
steps
(chosen randomly).\\

{\bf Observation 2}: The number of $k$-noncrossing tangled-diagrams
over $[n]$, (allowing isolated points), equals the number of simple
lattice walks in $x_1 \geq x_2 \geq \dots \geq x_{k-1} \geq 0$, from
the origin back to the origin, taking $n$ days, where at each day
the walker can either feel lazy and stay in place, or make {\it one}
unit step in any (legal) direction, or else feel energetic and make
any {\it two} consecutive steps
(chosen randomly).\\

These follow easily from the consideration in \cite{Reidys:07vac},
and are left as amusing exercises for the readers. The paper is
organized as follows: in Section~\ref{S:efficient} we consider
enumeration and computation using the holonomic framework \cite{Z}.
In Section~\ref{S:k} we validate that the formula, proved in
Section~\ref{S:efficient} for $k=2,3,4$, holds for arbitrary $k$.

\section{Efficient enumeration}\label{S:efficient}

Let $t_{k}(n)$ and $\tilde{t}_{k}(n)$ denote the numbers of
$k$-noncrossing tangled-diagrams with and without isolated points,
respectively. Furthermore let $f_k(m)$ denote the number of
$k$-noncrossing matchings over $m$ vertices or equivalently be the
number of ways of walking $n$ steps in $x_1 \geq x_2 \geq \dots \geq
x_{k-1} \geq 0$, from the origin back to the origin. Then, as shown
in \cite{Reidys:07vac}, $\tilde{t}_{k}(n)$ and $t_{k}(n)$ are given
by:
\begin{equation}\label{E:rel}
\tilde{t}_{k}(n)=\sum_{i=0}^{n} {{n} \choose {i}} f_k(2n-i) \quad
\text{\rm and}\quad t_{k}(n)=\sum_{i=0}^{n} {{n} \choose {i}}
\tilde{t}_{k} (n-i) \quad.
\end{equation}
Grabiner and Magyar proved an explicit determinant formula,
\cite{GM} (see also \cite{CDDSY}, eq.~$9$) that expresses the
exponential generating function of $f_k(n)$, for fixed $k$, as a
$(k-1) \times (k-1)$ determinant
\begin{eqnarray}
\label{E:ww1} \sum_{n\ge 0} f_{k}(2n)\cdot\frac{x^{2n}}{(2n)!} & = &
\det[I_{i-j}(2x)-I_{i+j}(2x)]|_{i,j=1}^{k-1} \ ,
\end{eqnarray}
where $I_m(2x)$ is the hyperbolic Bessel function:
\begin{equation}\label{E:bessel}
I_m(2x)=\sum_{j=0}^{\infty} { {x^{m+2j}} \over {j!(m+j)!} } \quad .
\end{equation}
Recall that a formal power series $G(x)$ is D-finite if it satisfies a
linear differential equation with polynomial coefficients.
For any $m$ the hyperbolic Bessel functions are $D$-finite \cite{S},
which is also called $P$-finite in \cite{Z}.
By general considerations, that we omit here, it is easy to
establish a priori bounds for the order of the recurrence, and for
the degrees of its polynomial coefficients, any {\it empirically}
derived recurrence (using the  command {\tt listtorec} in the
Salvy-Zimmerman Maple package {\tt gfun}, that we adapted to our own
needs in our own package {\tt TANGLE}), is {\it ipso facto}
rigorous. We derived explicit recurrences for $k=2,3,4$, and they
can be found in the webpage of this article. Also, once recurrences
are found, they are very efficient in extending the counting
sequences. In the same page one can find the sequences for
$T_{k}(n)$ for $1 \leq n \leq 1000$, for $k=2,3,4$, and the
sequences for $1 \leq n \leq 50$ for $k=5,6$ (using a variant of the
Grabiner-Magyar formula implemented in our Maple package { \tt
TANGLE} ).

Once the existence of a recursion is established, we can, for
$k=2,3,4$, employ the Birkhoff-Tritzinsky theory \cite{T:wimp,wimp}
and non-rigorous ``series analysis'' due to Zinn-Justin
\cite{zinn,zinn1}. This allows us to safely conjecture that, for any
fixed $k$, we have the following asymptotic formula:
\begin{equation}\label{E:A}
t_{k}(n) \quad \sim  \quad c_k \, \cdot \, n^{-((k-1)^2+(k-1)/2)} \,
(4(k-1)^2+2(k-1)+1)^n
  \quad \text{\rm for some $c_k>0$.}
\end{equation}
In the next Section we shall prove eq.~(\ref{E:A}) for arbitrary
$k$.

\section{Asymptotics of tangled-diagrams for arbitrary $k$}\label{S:k}


In Lemma~\ref{L:func1} we relate the generating functions of
$k$-noncrossing tangled diagrams $T_k(z)=\sum_{n}t_k(n)z^n$ and
$k$-noncrossing matchings \cite{CDDSY}
$F_k(z)=\sum_{n}f_k(2n)\,z^{2n}$. The functional equation derived
will be instrumental to prove eq.~(\ref{E:A}) for arbitrary $k$. For
this purpose we shall employ Cauchy's integral formula: let $D$ be a
simply connected domain and let $C$ be a simple closed positively
oriented contour that lies in $D$. If $f$ is analytic inside $C$ and
on $C$, except at the points $z_1,z_2,\dots, z_n$ that the interior
of $C$, then we have Cauchy's integral formula
\begin{equation}\label{E:cauchy}
\int_C\, f(z)dz=2\pi i\sum_{k=1}^{n}Res[f,z_k] \ .
\end{equation}
In particular, if $f$ has a simple pole at $z_0$, then
$Res[f,z_{0}]= \lim\limits_{z\rightarrow z_{0}}(z-z_0)f(z)$.
\begin{lemma}\label{L:func1}
Let $k\in\mathbb{N}$, $k\ge 2$ and $|z|<2$. Then we have
\begin{equation}
T_{k}\left(\frac{z^2}{1+z+z^2}\right)= \frac{1+z+z^2}{z+2}\,
F_{k}(z)\, .
\end{equation}
\end{lemma}
\begin{proof}
The relation between the number of $k$-noncrossing tangled-diagrams,
$t_k(n)$ and $k$-noncrossing matchings, $f_k(2m)$ given in
eq.~(\ref{E:rel}) implies $ t_{k}(n)=\sum_{r,\ell}{n \choose r}{n-r
\choose \ell}f_{k}(2n-2r-\ell)$. Expressing the combinatorial terms
by contour integrals we obtain
\begin{eqnarray*}
{n\choose r} & = & \frac{1}{2\pi i}\,
                        \oint_{|u|=\alpha}(1+u)^{n}u^{-r-1}du \\
f_{k}(2n-2r-\ell) & = & \frac{1}{2\pi i}\,
            \oint_{|z|=\beta_3}F_k(z)z^{-(2n-2r-\ell)-1}dz \\
t_{k}(n)&=&\sum_{r,\ell}{n \choose r}{n-r \choose \ell}
                                                  f_{k}(2n-2r-\ell) \\
&=&\frac{1}{(2\pi
i)^3}\sum_{r,\ell}\oint_{\substack{|v|=\beta_1\\|z|=
\beta_2\\|u|=\beta_3}}(1+u)^nu^{-r-1}(1+v)^{n-r}v^{-\ell-1} \ \times \\
& &  \qquad\qquad\qquad\qquad\qquad F_{k}(z)\,z^{-(2n-2r-\ell)-1}dv
\,du \, dz,
\end{eqnarray*}
where $\alpha,\beta_1,\beta_2,\beta_3$ are arbitrary small positive
numbers.
  Due to absolute convergence of the series we derive
\begin{eqnarray*}
t_{k}(n)&=&\frac{1}{(2\pi
i)^3}\sum_{r}\oint_{\substack{|v|=\beta_1\\|z|=
\beta_2\\|u|=\beta_3}}(1+u)^nu^{-r-1}
F_{k}(z)\,z^{-2n+2r-1}(1+v)^{n-r}v^{-1} \, \times \\
& & \qquad\qquad\qquad\qquad\qquad\qquad\qquad
\qquad\qquad\qquad\qquad
\sum_\ell\left (\frac{z}{v}\right)^{\ell}dv \,du \,dz,\\
\end{eqnarray*}
which is equivalent to
\begin{eqnarray*}
t_{k}(n)&=&\frac{1}{(2\pi i)^3}\sum_{r} \oint_{\substack{|u|=\beta_3\\
|z|=\beta_2}} (1+u)^nu^{-r-1}F_{k}(z)
\, z^{-2n+2r-1} \, \times \\
& & \qquad\qquad\qquad\qquad\qquad\qquad\qquad\qquad\qquad
\left(\oint_{|v|=\beta_1}\frac{(1+v)^{n-r}}{v-z}dv\right) du\,dz \ .
\label{E:int1}\\
\end{eqnarray*}
Since $v=z$ is the only (simple) pole in the integration domain,
eq.~(\ref{E:cauchy}) implies
$$
\oint_{|v|=\beta_1}\frac{(1+v)^{n-r}}{v-z}dv=2\pi i\,(1+z)^{n-r} \ .
$$
We accordingly obtain
\begin{equation}
t_{k}(n)=\frac{1}{(2\pi i)^2}\sum_{r} \oint_{\substack{|u|=\beta_3\\
|z|=\beta_2}} (1+u)^nu^{-r-1}F_{k}
(z)\,z^{-2n+2r-1}(1+z)^{n-r}du\,dz.
\end{equation}
Proceeding analogously w.r.t.~the summation over $r$ yields
\begin{eqnarray*}
t_{k}(n)&=&\frac{1}{(2\pi i)^2} \oint_{\substack{|u|=\beta_3\\
|z|=\beta_2}} (1+u)^nF_{k}
(z)\,z^{-2n-1}(1+z)^{n}u^{-1}\sum_{r}\frac{z^{2r}}
{u^{r}(1+z)^{r}}du \, dz\\
&=&\frac{1}{(2\pi i)^2} \oint_{|z|=\beta_2} F_{k}
(z)\,z^{-2n-1}(1+z)^{n}
\left(\oint_{|u|=\beta_3}(1+u)^n\frac{1}{u-\frac{z^{2}}
{1+z}}du\right)dz \ .
\end{eqnarray*}
Since $u=\frac{z^{2}}{1+z}$ is the only pole in the integration
domain, Cauchy's integral formula implies
$\oint_{|u|=\beta_3}(1+u)^n\frac{1}{u-\frac{z^{2}}{1+z}}du=2\pi
i\,(1+\frac{z^{2}}{1+z})^{n}$. We finally compute
\begin{eqnarray*}
t_{k}(n)&=&\frac{1}{2\pi i}\oint_{|z|=\beta_2} F_{k} (z)\,
z^{-1}z^{-2n}(1+z)^{n}
(1+\frac{z^{2}}{1+z})^{n}dz\\
&=&\frac{1}{2\pi i}\oint_{|z|=\beta_2} F_{k} (z)\,
z^{-1}\left(\frac{1+z+z^2}
{z^2}\right)^{n}dz \\
&=&\frac{1}{2\pi i}\oint_{|z|=\beta_2}
\frac{1+z+z^2}{z+2}F_{k}(z)\left(\frac
{z^2}{1+z+z^2}\right)^{-n-1}d{\left(\frac{z^2} {1+z+z^2}\right)}
\end{eqnarray*}
and the lemma follows from Cauchy's integral formula
\begin{equation}
T_{k}\left(\frac{z^2}{1+z+z^2}\right)=\frac{1+z+z^2} {z+2}F_{k}(z) \
.
\end{equation}
\end{proof}

\begin{theorem}
For arbitrary $k\in\mathbb{N}$, $k\ge 2$ the number of
tangled-diagrams is asymptotically given by
\begin{equation}
t_k(n) \sim c_k\, n^{-((k-1)^2 +\frac{k-1}{2})}\,
\left(4(k-1)^2+2(k-1)+1\right)^n \quad \text{\it where} \ c_k>0 \ .
\end{equation}
\end{theorem}
\begin{proof}
According to \cite{S,Z}, $F_k(x)=\sum_nf_k(2n)\,x^{2n}$ and $T_k(x)$
are both D-finite. Therefore both have a respective singular
expansion \cite{Flajolet}. We consider the following asymptotic
formula for $f_k(2n)$ \cite{Reidys:08k}: for arbitrary $k\ge 2$
\begin{equation}\label{E:k}
f_k(2n)\sim n^{-((k-1)^2 +\frac{k-1}{2})}\, \left( 2(k-1)\right)^{2n} \
.
\end{equation}
Eq.~(\ref{E:k}) allows us to make two observations. First $F_k(x)$
has the positive, real, dominant singularity, $\rho_k=(2(k-1))^{-1}$
and secondly, in view of the subexponential factor $n^{-((k-1)^2
+\frac{k-1}{2})}$:
\begin{eqnarray}\label{E:e1}
F_k(z)  & = &  O\left(\left( z -\rho_k\right)^{((k-1)^2
+\frac{k-1}{2})-1} \right),\qquad \text{\rm as $z\rightarrow
\rho_k$.}
\end{eqnarray}
According to Lemma~\ref{L:func1} we have
\begin{equation}\label{E:oha}
T_k\left(\frac{z^2}{z^2+z+1}\right) \, = \, \frac{z^2+z+1} {z+2}\,
F_{k}(z) \ ,
\end{equation}
where $\vert z\vert \le \rho_k\le \frac{1}{2}+\epsilon$,
$\epsilon>0$ is arbitrarily small and the function $\vartheta(z)=
\frac{z^2}{z^2+z+1}$ is regular at $z=\rho_k$. Since the composition
$H(\eta(z))$ of a D-finite function $H$ and a rational function
$\eta$, where $\eta(0)=0$ is D-finite \cite{S}, the functions
$T_k(\vartheta(z))$ and $F_k(z)$ have singular expansions.
Eq.~(\ref{E:oha}) and eq.~(\ref{E:k}) imply using Bender's method
($F_k(z)$ satisfies the ``ratio test'')
\cite{Bender}
\begin{equation}\label{E:wichtig}
[z^n]\, T_k(\vartheta(z))\,\sim  \,
\frac{\rho_k^2+\rho_k+1}{\rho_k+2}\, [z^n]\, F_k(z)\,\sim\,
\frac{\rho_k^2+\rho_k+1}{\rho_k+2}\; n^{-((k-1)^2 +\frac{k-1}{2})}\,
\left( \rho_k^{-1}\right)^{2n} \ .
\end{equation}
Eq.~(\ref{E:wichtig}) implies that
$\tau_k=\frac{\rho_k^2}{\rho_k^2+\rho_k+1}$ is the positive, real,
dominant singularity of $T_k(z)$. Indeed, Pringsheim's Theorem
\cite{Titmarsch:39} guarantees the existence of a positive, real,
dominant singularity of $T_k(z)$, denoted by $\tau_k$. For $0\le
x\le 1$ the mapping $x\mapsto \frac{x^2}{x^2+x+1}$ is strictly
increasing and continuous, whence
$\tau_k=\frac{\zeta^2}{\zeta^2+\zeta+1}$ for some $0< \zeta \le 1$.
In view of eq.~(\ref{E:oha}), $\zeta$ is the dominant, positive,
real singularity of $F_k(z)$, i.e.~$\zeta=\rho_k$. Accordingly,
\begin{equation}
[z^n]\,T_k(z) \sim \theta(n) \,
\left(\frac{\rho_k^2}{\rho_k^2+\rho_k+1}\right)^n \ .
\end{equation}
We shall proceed by analyzing $T_k(z)$ at dominant singularities. We
observe that any dominant singularity $v$ can be written as
$v=\vartheta(\zeta)$. Let $S_{T_k}(z-\vartheta(\zeta))$ denote the
singular expansion of $T_k(z)$ at $v=\vartheta(\zeta)$. Since
$\vartheta(z)$ is regular at $\zeta$, $T_k(\vartheta(z))$ we have
the {\it supercritical} case of singularity analysis
\cite{Flajolet}: given $\psi(\phi(z))$, $\phi$ being regular at the
singularity of $\psi$, the singularity-type of the composition is
that of $\psi$. Indeed, we have
\begin{eqnarray*}\label{E:DD}
T_k(\vartheta(z)) & = & O(S_{T_k}(\vartheta(z)-\vartheta(\zeta)))
            \qquad \text{\rm as $\vartheta(z)\rightarrow \vartheta(\zeta)$.}\\
                  & = & O(S_{T_k}(z-\zeta)) \ \ \, \qquad \qquad \text{\rm as
$z\rightarrow \zeta$}.
\end{eqnarray*}
Eq.~(\ref{E:oha}) provides the following interpretation for
$T_k(\vartheta(z))$ at $z=\zeta$:
\begin{eqnarray*}
T_k(\vartheta(z)) = O(F_k(z)) \qquad \qquad \text{\rm as
$z\rightarrow \zeta$},
\end{eqnarray*}
from which we can conclude that $T_k(z)$ has at $v=\vartheta(\zeta)$
exactly the same subexponential factors as $F_k(z)$ at $\zeta$. We
next prove that $\tau_k$ is the unique dominant singularity of
$T_k(z)$. Suppose $v=\vartheta(\zeta)$ is an additional dominant
singularity of $T_k(z)$. The key observations is
\begin{equation}\label{E:rule}
\forall \;\zeta\in\mathbb{C}\setminus \mathbb{R};\qquad
\vartheta(\zeta) = \tau_k \quad \Longrightarrow \quad \vert
\zeta\vert < \rho_k \ .
\end{equation}
Eq.~(\ref{E:rule}) implies that if $v$ exists then $\zeta$ is a
singularity of $F_k(z)$ of modulus strictly smaller than $\rho_k$,
which is impossible. Therefore $\tau_k$ is unique and we derive
\begin{equation}
[z^n]\,T_k(z) \sim c_k \, n^{-((k-1)^2 +\frac{k-1}{2})}\,
\left(\frac{\rho_k^2}{\rho_k^2+\rho_k+1}\right)^n \quad \text{\rm
for some $c_k>0$}
\end{equation}
and the theorem follows.
\end{proof}
{\bf Acknowledgments.}
We are grateful to Emma Y.~Jin for helpful discussions. This work
was supported by the 973 Project, the PCSIRT Project of the Ministry
of Education, the Ministry of Science and Technology, and the
National Science Foundation of China. The fourth author is supported
in part by the USA National Science Foundation.
\bibliographystyle{amsplain}

\end{document}